\documentclass[reqno,11pt]{article}
\usepackage{amssymb}
\usepackage{amsmath}
\usepackage{amsthm}
\usepackage{geometry} 

\newtheorem{Theorem}{Theorem } [section]
\newtheorem{lemma}[Theorem]{Lemma}

\newtheorem{remark}{Remark}

\numberwithin{equation}{section}

\newcommand{\vertiii}[1]{{\vert\kern-0.25ex\vert\kern-0.25ex\vert #1 
    \vert\kern-0.25ex\vert\kern-0.25ex\vert}}

\DeclareMathOperator{\supp}{supp}

\date{}  

\title{Dispersive blow-up for the fifth order Korteweg-de Vries equation on the line}
\author{Eddye Bustamante, Jos\'e Jim\'enez Urrea and Jorge Mej\'{\i}a}

\author{
    Eddye Bustamante \\
    \textit{Universidad Nacional de Colombia, Sede Medell\'in} \\
    \texttt{eabusta0@unal.edu.co} 
    \and
    José Jiménez Urrea \\
    \textit{Universidad Nacional de Colombia, Sede Medell\'in} \\
    \texttt{jmjimene@unal.edu.co} 
    \and
    Jorge Mej\'ia \\
    \textit{Universidad Nacional de Colombia, Sede Medell\'in} \\
    \texttt{jemejia@unal.edu.co} 
}

\begin{document}

%
%

\maketitle

\begin{abstract} In this work we establish a dispersive blow-up result for the initial value problem (IVP) for the fifth order Korteweg-de Vries equation
\begin{align*}
\left. \begin{array}{rlr}
u_t+\partial_x^5 u+u\partial_x u&\hspace{-2mm}=0,&\quad x\in\mathbb R,\; t>0,\\
u(x,0)&\hspace{-2mm}=u_0(x),&
\end{array} \right\}
\end{align*}
To achieve this, we prove a local well-posedness result in Bourgain spaces of the type $X^{s,b}$ for appropriate values of $s$ and $b$, along with a regularity property for the nonlinear part of that solution. This property enables the construction of initial data that leads to the dispersive blow-up phenomenon.
\end{abstract}

\setcounter{page}{1}
\pagenumbering{arabic}

\section{Introduction}

\footnotetext[1]{\textbf{2020 Mathematics Subject Classification:} 35Q53}

In this article we consider the initial value problem (IVP) corresponding to the fifth order Korteweg-de Vries (KdV) equation 
\begin{align}
\left. \begin{array}{rlr}
u_t+\partial_x^5 u+u\partial_x u&\hspace{-2mm}=0,&\quad x\in\mathbb R,\; t>0,\\
u(x,0)&\hspace{-2mm}=u_0(x).&
\end{array} \right\}\label{maineq}
\end{align}

The well-posedness of the Cauchy problem associated with fifth-order KdV-type equations in the functional context of Sobolev spaces has been addressed by various authors in recent years (see, for example, \cite{BBM1972}, \cite{BJM2015}, \cite{GKK2013}, \cite{KP2015}, \cite{K2008a}, \cite{K2008b}, \cite{L1995}, and references therein).\\

The IVP \eqref{maineq} is a particular case of the IVP, studied by Kato in \cite{K2011} and \cite{K2013}, for Kawahara equation
\begin{align}
\left. \begin{array}{rlr}
u_t+\alpha \partial_x^5 u+ \beta \partial_x^3 u + \gamma\partial_x (u^2)&\hspace{-2mm}=0,&\quad x\in\mathbb R,\; t>0,\\
u(x,0)&\hspace{-2mm}=u_0(x),&
\end{array} \right\}\label{kawahara}
\end{align}
where $\alpha$, $\beta$, and $\gamma$ are real constants with $\alpha\neq 0$, $\gamma\neq 0$.\\

The equation in \eqref{kawahara} models the capillary waves on a shallow layer and the magneto-sound waves in plasma physics (see \cite{K1972}).\\

In \cite{K2011} Kato proves, in the context of Bourgain spaces, that the IVP \eqref{kawahara} is locally well-posed (LWP) with initial data in the Sobolev space $H^s(\mathbb R)$ with $s\geq-2$. Moreover, he proves ill-posedness for $s<-2$.\\

In \cite{K2013}, using a modified Bourgain space and applying the I-method, Kato manages to show that the IVP \eqref{kawahara} is globally well-posed in $H^s(\mathbb R)$ with $s\geq-\frac{38}{21}$.\\

Since in our work we intend to prove a regularity property of the nonlinear part of the solution of the IVP \eqref{maineq}, using an extra gain of regularity of the bilinear form $\partial_x(uv)$, measured in an appropriate Bourgain space $X^{s,b}$, then our first result is to establish again the local well-posedness of the IVP \eqref{maineq} in $X^{s,b}$, when the initial datum $u_0$ belongs to $H^s(\mathbb R)$ with $s\geq 0$.\\

Let us remember that, formally, $u$ is a solution of the IVP \eqref{maineq} if and only if
\begin{align}
u(t) = W(t)u_0 + \int_0^t W(t-t') (-\frac12 \partial_x (u(t'))^2) dt', \label{inteq}
\end{align}

where $\{W(t)\}_{t\in \mathbb R}$ is the group of unitary operators in the Sobolev spaces $H^s\equiv H^s(\mathbb R)$ associated to the linear part of the equation in \eqref{maineq}; i.e.,
\begin{align}
[W(t)u_0]^\wedge(\xi) = e^{-it\xi^5} \widehat  u_0(\xi). \label{group}
\end{align}
(Here, $u(t)\equiv u(\cdot_x,t)$ and $\widehat u_0$ is the Fourier transform of the space function $u_0$).\\

We are going to work in the context of the Bourgain spaces $X^{s,b}\equiv X^{s,b}(\mathbb R^2)$ of tempered distributions $u$ in $\mathbb R^2_{xt}$ such that
\begin{align}
\|\langle \xi \rangle^s \langle \tau + \xi^5 \rangle^b \widehat  u (\xi,\tau)\|_{L^2_{\xi\tau}} < \infty,\label{Bnorm}
\end{align}

where $s\in\mathbb R$, and $b=\frac12^-$; i.e., there exists $\epsilon>0$ such that $\frac12-\epsilon<b<\frac12$. Here $\langle\cdot\rangle$ is an abbreviation for $1+|\cdot|$, $\widehat u$ denotes the Fourier transform of the distribution $u$, and $\xi$ and $\tau$ are the variables in the frequency space corresponding to $x$ and $t$, respectively.\\


We prove that, for suitable values of $s$ and $b$, and for $u_0\in H^s(\mathbb R)$, there exists $u\in X^{s,b}(\mathbb R^2)$, continuous from $\mathbb R_t$ with values in $H^s(\mathbb R_x)$, i.e., $u\in X^{s,b}(\mathbb R^2) \cap C(\mathbb R_t; H^s(\mathbb R_x))$, such that, for some $T>0$ and $t\in[0,T]$, equation \eqref{inteq} holds. In order to achieve this, we introduce a bump function $\eta$ of the time variable $t$ such that $\eta(\cdot_t)\in C_0^\infty(\mathbb R_t)$, $\supp \eta\subset [-1,1]$, and $\eta\equiv 1$ in $[-\frac12,\frac12]$, and we consider the integral equation

\begin{align}
u(t) = \eta(t) W(t)u_0 + \eta(t) \int_0^t W(t-t') F_T(u(t')) dt', \label{inteqtr}
\end{align}

where $F_T(u(t')) = \eta(\frac{t'}{2T}) (-\frac12 \partial_x (u(t'))^2)$. In a precise manner, we prove the following result.

\begin{Theorem}\label{localsolintime} (Local well posedness of the IVP \eqref{maineq}). Let $s\geq 0$, and $u_0\in H^s(\mathbb R)$. For $\frac25<b<\frac12$, there exist $T\in(0,\frac12]$, and a unique $u\in X^{s,b}\cap C(\mathbb R_t; H^s(\mathbb R_x))$ solution of \eqref{inteqtr}. (The restriction of $u$ to $\mathbb R\times[0,T]$ is a local solution in time of the IVP \eqref{maineq}).
\end{Theorem}

It can be seen that the result of Theorem \ref{localsolintime} is also global in time. \\

Our second result is related to the extra regularity of the nonlinear part of the solution of the IVP \eqref{maineq}. This type of result was studied in the case of the IVP associated to the generalized Korteweg-de Vries (KdV) equation
\begin{align}
\left. \begin{array}{rlr}
u_t+\partial_x^3 u + u^k \partial_x u&\hspace{-2mm}=0,&\quad k\in \mathbb Z^+,\quad k\geq2,\quad x\in\mathbb R,\; t>0,\\
u(x,0)&\hspace{-2mm}=u_0(x),&
\end{array} \right\}\label{kdv}
\end{align}

by Linares and Scialom in \cite{LS1993}. Specifically, if we write the solution $u(x,t)$ of the IVP \eqref{kdv} in the form
\begin{align*}
u(x,t) & = [V(t)u_0](x) + \left[ \int_0^t V(t-t')[-u^k(t')\partial_x u(t')] dt' \right] (x)\equiv u_1(x,t) + u_2(x,t),
\end{align*}

where $\{V(t)\}_{t\in\mathbb R}$ is the group associated to the linear part of the equation in \eqref{kdv}, Linares and Scialom prove that the integral part $u_2$ of the solution $u$, is smoother than $u$, when the initial data $u_0$ is taken in $H^s(\mathbb R)$ with $s\geq 1$. Precisely they demonstrate that if $u_0\in H^s(\mathbb R)$ with $s\geq 1$, then $u_2\in C(\mathbb R_t; H^{s+1}(\mathbb R_x))$ and for all $T>0$ $\sup_x \|\partial_x ^{s+2} u_2\|_{L^2_T}<+\infty$. This result was obtained, using local smoothing effects of Kato type.\\

Later, in \cite{LPS2017}, Linares, Ponce and Smith construct an initial data $u_0\in H^{\frac 32^-}(\mathbb R)$ such that $|x|^\beta u_0 \in L^2(\mathbb R)$ for $\beta\in(0,\frac34)$ in such a way that for the IVP associated to the Korteweg-de Vries equation ($k=1$ in the equation in \eqref{kdv}) the nonlinear part of the solution of the IVP, given by
\begin{align*}
u_2(x,t):= \left[ \int_0^t V(t-t') (-u(t')\partial_x u(t') dt') \right] (x)
\end{align*}

is such that $u_2\in C([0,T]; H^{\frac32^- +\frac16}(\mathbb R))$. In other words, $u_2(t)$ is smoother than $u_0$ for every $t\in[0,T]$. In addition to using regularizing effects, in \cite{LPS2017} the authors make use of the local well-posedness of the IVP in weighted Sobolev spaces $H^{\frac32-\epsilon}(\mathbb R) \cap L^2(|x|^{\frac34-2\epsilon} dx)$.\\

In the case of systems of differential equations, recently, Linares and Palacios in \cite{LP2019} prove that for a Schrödinger-Korteweg-de Vries system the nonlinear part of the solution of the IVP associated to this system is smoother than the solution. As in \cite{LPS2017}, they use smoothing properties, and persistence properties for the solutions of the IVP in weighted Sobolev spaces.\\

Using Bourgain spaces, in \cite{ET2016}, Erdoğan and Tzirakis prove that the nonlinear part of the solution of an initial boundary value problem for the cubic nonlinear Schrödinger equation is smoother than the initial data. This result is a consequence of some gain of regularity of the nonlinear term $|u|^2 u$, measured in the norm of certain Bourgain space. With the help of this approach, in the case of the fifth KdV equation, we manage to prove that the nonlinear part of the solution of the IVP \eqref{maineq} is smoother than the initial data $u_0$. More precisely, we obtain the following result
\begin{Theorem}\label{more_reg} (Regularity gain of the nonlinear part of the Solution). Let $s\geq 0$, $u_0\in H^s(\mathbb R)$, $\frac25<b<\frac12$, and $0\leq a\leq 10b-4$. On the other hand, let $T\in (0,\frac 12]$ and $u\in X^{s,b}\cap C(\mathbb R_t;H^s(\mathbb R_x))$ be as in Theorem \ref{localsolintime}. Then
\begin{align}
\left\| \int_0^{\cdot_t} W(\cdot_t{\scriptstyle -}t')(-\frac12 \partial_x(u(t'))^2) dt'\right\|_{C([0,T];H^{s+a}(\mathbb R))} \leq C \|u\|^2_{X^{s,b}}.\label{more_reg_eq1}
\end{align}
\end{Theorem}

Our main objective in this paper is to establish dispersive blow-up for the solutions of the IVP \eqref{maineq}. The pioneering work in this direction for dispersive equations is due to Benjamin, Bona, and Mahony in \cite{BBM1972}. In that paper, the authors show that there exist smooth initial data $u_0$, for which the solution of the linearized KdV equation associated with these initial data can develop point singularities in finite time. Other results related to dispersive blow-up for different nonlinear dispersive models can be found in \cite{BPSS2014}, \cite{BS2010}, \cite{MP2023}, and the references therein.\\

In \cite{BS1993}, Bona and Saut prove the existence of dispersive blow-up for the generalized KdV equation. The idea behind the Bona and Saut proof is to observe that the nonlinear part of the solution is smoother than the solution, which allows to focus attention on the construction of initial data that gives rise to the formation of point singularities for the solution of the linearized IVP. Following this idea, in \cite{LPS2017} and \cite{LP2019} it is possible to establish the dispersive blow-up in the case of the generalized KdV equation and in the case of a Schrödinger-KdV system, respectively.\\

Finally, using the regularity gain in the Theorem \ref{more_reg} and by imitating the procedure followed in \cite{LPS2017} (see proof of Theorem 1.4 in \cite{LPS2017}) we prove a result on dispersive blow-up, which is stated as follows.
\begin{Theorem}\label{blowup} (Dispersive blow-up). There exists $u_0\in H^1(\mathbb R)\cap C^{\infty}(\mathbb R)$ with $\|u_0\|_{H^1}\ll 1$ so that the solution $u(\cdot,t)$ of the IVP \eqref{maineq} satisfies the following properties:
\begin{enumerate}
\item[(i)] $u(\cdot,t)\in C^1(\mathbb R)$ for $t\in(0,+\infty) \setminus \mathbb Z^+$;
\item[(ii)] $u(\cdot,t)\in C^1(\mathbb R\setminus \{0\}) \setminus C^1(\mathbb R)$ for $t\in \mathbb Z^+$.
\end{enumerate}

\end{Theorem}

\section{Linear and non linear estimates}

In this section we establish some estimates for the norms in the spaces $X^{s,b}$ of the terms on the right hand side of the integral equation \eqref{inteqtr}. Throughout this section the function $\eta(\cdot_t)$ is as before.

\begin{lemma}\label{le1} Let $s\in\mathbb R$, $b\in\mathbb R$. Then there exists $C>0$ such that, for every function $u_0\in H^s(\mathbb R_x)$,
\begin{align}
\|\eta(\cdot_t)[W (\cdot_t)u_0](\cdot_x)\|_{X^{s,b}(\mathbb R^2)}\leq C\|u_0\|_{H^s(\mathbb R_x)}.\label{le1_eq1}
\end{align}
\end{lemma}

\begin{proof} Let us observe that
\begin{align*}
\left(\eta(\cdot_t)[W (\cdot_t)u_0](\cdot_x)\right)^\wedge(\xi,\tau)=C\int_{\mathbb R_t} e^{-i(\tau+\xi^5)t}\eta(t)\widehat u_0(\xi) dt=C\widehat u_0(\xi)\widehat \eta(\tau+\xi^5).
\end{align*}

Therefore, taking into account that $\widehat \eta\in S(\mathbb R)$,
\begin{align*}
\|\eta(\cdot_t)[W (\cdot_t)u_0](\cdot_x)\|_{X^{s,b}(\mathbb R^2)}&=C\left( \iint_{\mathbb R^2}\langle\xi\rangle^{2s}\langle\tau+\xi^5\rangle^{2b}|\widehat u_0(\xi)|^2 |\widehat\eta(\tau+\xi^5)|^2d\tau d\xi\right)^{\frac12}\\
&=C\left( \int_{\mathbb R_{\tau}} \langle\tau\rangle^{2b}|\widehat\eta(\tau)|^2\right)^{\frac12}\|u_0\|_{H^s(\mathbb R_x)}=C\|u_0\|_{H^s(\mathbb R_x)}.
\end{align*}
\end{proof}

\begin{remark} For $s,b\in\mathbb R$, let us denote by $H^{s,b}$ the anisotropic Sobolev space, defined by
$$H^{s,b}:=\{f(\cdot_x,\cdot_t)\in S'(\mathbb R^2):\|f\|_{H^{s,b}}:=\|\langle\xi\rangle^s\langle\tau\rangle^b\widehat f(\xi,\tau)\|_{L^2_{\xi \tau}}<\infty\}.$$

It can be seen that
\begin{align}
\|f(\cdot_x,\cdot_t)\|_{X^{s,b}}=\|W ({\scriptstyle-}\cdot_t)f(\cdot_x,\cdot_t)\|_{H^{s,b}}.\label{r1_eq1}
\end{align}
\end{remark}

\begin{remark} In \cite{GTV1997}, Ginibre, Tsutsumi, and Velo proved that if $L$ is the operator defined by
$$(Lf)(t):=\eta\left(\tfrac tT\right)\int_0^tf(t')dt',$$
for fixed $T\in(0,1]$, then, for $-\frac12<b'\leq 0\leq b\leq b'+1$, it is true that
\begin{align}
\|Lf\|_{H^b(\mathbb R_t)}\leq C T^{1-b+b'}\|f\|_{H^{b'}(\mathbb R_t)},\label{r2_eq1}
\end{align}
where $C>0$ is independent of $f$.
\end{remark}

\begin{lemma}\label{le2}
For $s\in\mathbb R$, $-\frac12<b'\leq 0\leq b\leq b'+1$, $0<T\leq 1$, there exists $C>0$, such that for every $f\in X^{s,b'}$,
\begin{align}
\left\|\eta\left(\tfrac{\cdot_t}T \right)\int_0^{\cdot_t} [W (\cdot_t{\scriptstyle -}t')f(t')](\cdot_x)dt'\right\|_{X^{s,b}}\leq C T^{1-b+b'}\|f\|_{X^{s,b'}}\label{le2_eq1}
\end{align}
\end{lemma}

\begin{proof} Using \eqref{r1_eq1} and \eqref{r2_eq1} we obtain
\begin{align}
\notag &\left\|\eta\left(\tfrac{\cdot_t}T \right)\int_0^{\cdot_t} [W (\cdot_t{\scriptstyle-}t')f(t')](\cdot_x)dt'\right\|_{X^{s,b}}\\
\notag&=\left\| W ({\scriptstyle-}\cdot_t)\left\{ \eta\left(\tfrac{\cdot_t}T \right)\int_0^{\cdot_t} [W (\cdot_t{\scriptstyle-}t')f(t')](\cdot_x)dt'\right\}\right\|_{H^{s,b}}\\
\notag&=\left(\int_{\mathbb R_\xi}\langle\xi\rangle^{2s}\int_{\mathbb R_\tau}\langle\tau\rangle^{2b} \left|\left\{  \eta\left(\tfrac{\cdot_t}T \right)\left(\int_0^{\cdot_t} [W (-t')f(t')](\cdot_x)dt'\right)^{\wedge_x}\right\}^{\wedge_t}(\xi,\tau)\right|^2d\tau d\xi\right)^{\frac12}\\
\notag&=\left(\int_{\mathbb R_\xi}\langle\xi\rangle^{2s} \left\| L\left(\left( [W ({\scriptstyle-}\cdot_t)f(\cdot_t)](\cdot_x)\right)^{\wedge_x}(\xi)\right)\right\|^2_{H^b(\mathbb R_t)} d\xi\right)^{\frac12}\\
&\leq CT^{1-b+b'}\left(\int_{\mathbb R_\xi}\langle\xi\rangle^{2s}\left\|\left([W ({\scriptstyle-}\cdot_t)f(\cdot_t)](\cdot_x) \right)^{\wedge_x}(\xi) \right\|^2_{H^{b'}(\mathbb R_t)}d\xi \right)^{\frac12}.\label{le2_eq2}
\end{align}

Taking into account that
\begin{align*}
\left\|\left([W ({\scriptstyle -}\cdot_t)f(\cdot_t)](\cdot_x) \right)^{\wedge_x}(\xi) \right\|^2_{H^{b'}(\mathbb R_t)}=\int_{\mathbb R_\tau}\langle\tau\rangle^{2b'}\left|\left(e^{it\xi^5}[f(\cdot_x,t)]^{\wedge_x}(\xi) \right)^{\wedge_t}(\tau)\right|^2d\tau,
\end{align*}
and since

\begin{align*}
\left(e^{it\xi^5}[f(\cdot_x,t)]^{\wedge_x}(\xi) \right)^{\wedge_t}(\tau)&=C\int_{\mathbb R_t}e^{-it(\tau-\xi^5)}\int_{\mathbb R_x}e^{-i\xi x}f(x,t)dxdt,
\end{align*}
we conclude that
\begin{align*}
&\left\|\left([W ({\scriptstyle -}\cdot_t)f(\cdot_t)](\cdot_x) \right)^{\wedge_x}(\xi) \right\|^2_{H^{b'}(\mathbb R_t)}\\
&=C\int_{\mathbb R_\tau}\langle\tau\rangle^{2b'}|\widehat f(\xi,\tau-\xi^5)|^2d\tau=C\int_{\mathbb R_\tau}\langle\tau+\xi^5\rangle^{2b'}|\widehat f(\xi,\tau)|^2d\tau.
\end{align*}

From \eqref{le2_eq2}, we obtain
\begin{align*}
\left\|\eta\left(\tfrac{\cdot_t}T \right)\int_0^{\cdot_t} [W (\cdot_t{\scriptstyle -}t')f(t')](\cdot_x)dt'\right\|_{X^{s,b}}\leq CT^{1-b+b'}\|f\|_{X^{s,b'}}.
\end{align*}
\end{proof}

\begin{lemma}\label{le3} Let $s\in\mathbb R$, $b_1,b_2$ such that $-\frac12<b_1<b_2<\frac12$, and $T\in(0,1]$. Then there exists $C>0$ such that, for every $F\in X^{s,b_2}$,
\begin{align}
\left\|\eta\left(\tfrac{\cdot_t}T\right)F(\cdot_x,\cdot_t) \right\|_{X^{s,b_1}}\leq CT^{b_2-b_1}\|F\|_{X^{s,b_2}}.\label{le3_eq1}
\end{align}
\end{lemma}

\begin{proof} Let us prove \eqref{le3_eq1} for $b_1=0$; i.e.,
\begin{align}
\left\|\eta\left(\tfrac{\cdot_t}T\right)F(\cdot_x,\cdot_t) \right\|_{X^{s,0}}\leq CT^{b_2}\|F\|_{X^{s,b_2}}.\label{le3_eq2}
\end{align}

The general case follows from \eqref{le3_eq2}, using an interpolation argument. Estimative \eqref{le3_eq2} is consequence of the following statements:
\begin{enumerate}
\item[(i)] $\left\|\eta\left(\frac{\cdot_t}T \right)F(\cdot_x,\cdot_t) \right\|_{X^{s,0}}\leq CT^{b_2}\left\| \eta\left(\frac{\cdot_t}T\right)F(\cdot_x,\cdot_t)\right\|_{X^{s,b_2}}$,
\item[(ii)] $\left\|\eta\left(\frac{\cdot_t}T \right)F(\cdot_x,\cdot_t) \right\|_{X^{s,b_2}}\leq C\|F\|_{X^{s,b_2}}$;
\end{enumerate}
where $C>0$ is independent of $T$.\\

Let us prove (i):
\begin{align*}
\left\|\eta\left(\tfrac{\cdot_t}T \right)F(\cdot_x,\cdot_t) \right\|^2_{X^{s,0}}&=\int_{\mathbb R}\int_{-\infty}^{+\infty}\langle\xi\rangle^{2s}\left|[\eta\left( \tfrac{\cdot_t}T\right)F]^{\wedge}(\xi,\tau) \right|^2d\tau d\xi\\
&=\int_{\mathbb R}\langle\xi\rangle^{2s}\int_{-\infty}^{+\infty} \left|\eta\left( \tfrac tT\right)[F(\cdot,t)]^{\wedge_x}(\xi) \right|^2dt d\xi\\
&=\int_{\mathbb R_t} \left|\eta\left( \tfrac tT\right)\right|^2\int_{\mathbb R_\xi}\left| e^{it\xi^5} \langle\xi\rangle^{s} [F(\cdot,t)]^{\wedge_x}(\xi) \right|^2d\xi dt\\
&=C\int_{-T}^T \left|\eta\left( \tfrac tT\right)\right|^2\int_{\mathbb R_x}\left|[W (-t)J^s(F(\cdot,t))] (x)\right|^2dx dt
\end{align*}
where
$$[J^s(F(\cdot,t))]^{\wedge_x}(\xi):=\langle\xi\rangle^s[F(\cdot,t)]^{\wedge_x}(\xi).$$

Hence
\begin{align*}
\left\|\eta\left(\tfrac{\cdot_t}T \right)F(\cdot_x,\cdot_t) \right\|^2_{X^{s,0}}=C\int_{-T}^T\int_{\mathbb R_x}\left| [W (-t)J^s(\eta\left(\tfrac tT\right)F(\cdot,t))](x)\right|^2dxdt.
\end{align*}

Since $0<2b_2<1$, $\frac1{2b_2}>1$. Defining $p:=\frac1{2b_2}$, and $q$ as the conjugate exponent of $p$; i.e., $q:=\frac1{1-2b_2}$, by Hölder's inequality in the $t$ variable, we obtain
\begin{align*}
&\left\|\eta\left(\tfrac{\cdot_t}T \right)F(\cdot_x,\cdot_t) \right\|^2_{X^{s,0}}\\
&\leq C T^{2b_2}\left\{\int_{-\infty}^{+\infty}\left(\int_{\mathbb R}\left|[W (-t)J^s(\eta\left(\tfrac tT \right)F(\cdot,t))](x) \right|^2dx \right)^{\frac1{1-2b_2}} dt\right\}^{1-2b_2}\\
&\leq C T^{2b_2}\int_{\mathbb R}\left(\int_{-\infty}^{+\infty}\left|[W (-t)J^s(\eta\left(\tfrac tT \right)F(\cdot,t))](x) \right|^{\frac2{1-2b_2}}  dt\right)^{1-2b_2}dx\\
&=CT^{2b_2}\int_{\mathbb R} \left\|[W ({\scriptstyle -}\cdot_t)J^s(\eta\left(\tfrac {\cdot_t}T \right)F(\cdot,\cdot_t))](x) \right\|^2_{L^{\frac2{1-2b_2}}(\mathbb R_t)}dt.
\end{align*}

Taking into account that $H^{b_2}(\mathbb R)\hookrightarrow L^{\frac2{1-2b_2}}(\mathbb R)$ for $0\leq b_2<\frac12$, it follows that
\begin{align*}
&\left\|\eta\left(\tfrac{\cdot_t}T \right)F(\cdot_x,\cdot_t) \right\|^2_{X^{s,0}}\\
&\leq C T^{2b_2}\int_{\mathbb R}\int_{-\infty}^{+\infty}\langle\tau\rangle^{2b_2} \left|\left([W ({\scriptstyle-}\cdot_t)J^s(\eta\left(\tfrac {\cdot_t}T \right)F(\cdot,\cdot_t))](x)\right)^{\wedge_t} (\tau)\right|^2d\tau dx\\
&=CT^{2b_2}\int_{-\infty}^{+\infty}\langle\tau\rangle^{2b_2}\int_{\mathbb R_{\xi}} \left|\left\{\left([W ({\scriptstyle-}\cdot_t)J^s(\eta\left(\tfrac {\cdot_t}T \right)F(\cdot,\cdot_t))](\cdot_x)\right)^{\wedge_t}(\tau)\right\}^{\wedge_x}(\xi) \right|^2d\xi d\tau\\
&=CT^{2b_2}\int_{-\infty}^{+\infty}\langle\tau\rangle^{2b_2}\int_{\mathbb R_{\xi}} \left|\left(e^{it\xi^5}\langle\xi\rangle^s [\eta\left(\tfrac {\cdot_t}T \right)F(\cdot_x,\cdot_t)]^{\wedge_x}(\xi)\right)^{\wedge_t}(\tau) \right|^2d\xi d\tau\\
&=CT^{2b_2}\int_{-\infty}^{+\infty}\langle\tau\rangle^{2b_2}\int_{\mathbb R_{\xi}}\langle\xi\rangle^{2s} \left| [\eta\left(\tfrac {\cdot_t}T \right)F(\cdot_x,\cdot_t)]^{\wedge_{xt}}(\xi,\tau-\xi^5) \right|^2d\xi d\tau\\
&=CT^{2b_2}\int_{\mathbb R_\xi}\langle\xi\rangle^{2s}\int_{-\infty}^{+\infty}\langle\tau+\xi^5\rangle^{2b_2} \left| [\eta\left(\tfrac {\cdot_t}T \right)F(\cdot_x,\cdot_t)]^{\wedge_{xt}}(\xi,\tau) \right|^2d\tau d\xi\\
&=CT^{2b_2}\left\|\eta\left(\tfrac{\cdot_t}T \right)F(\cdot_x,\cdot_t) \right\|^2_{X^{s,b_2}},
\end{align*}
and estimative (i) is proved.\\

Now we prove (ii):
\begin{align*}
\left\|\eta\left(\tfrac{\cdot_t}T \right)F(\cdot_x,\cdot_t) \right\|^2_{X^{s,b_2}}=\int_{\mathbb R}\langle\xi\rangle^{2s}\int_{-\infty}^{+\infty}\langle\tau+\xi^5\rangle^{2b_2}\left| [\eta\left(\tfrac{\cdot_t}T \right)F(\cdot_x,\cdot_t)]^{\wedge}(\xi,\tau)\right|^2 d\tau d\xi.
\end{align*}

Let us observe that
\begin{align*}
&[\eta\left(\tfrac{\cdot_t}T\right)F(\cdot_x,\cdot_t)]^\wedge(\xi,\tau)\\
&=C\left[\eta\left(\tfrac {\cdot_t}T\right) (F(\cdot_x,t))^{\wedge_x}(\xi) \right]^{\wedge_t}(\tau)=C\left\{\left[\eta\left(\tfrac {\cdot_t}T\right)\right]^{\wedge_t}\ast_\tau \widehat F(\xi,\cdot_\tau)\right\}(\tau).
\end{align*}

But,
$$\left[\eta\left(\tfrac {\cdot_t}T\right)\right]^{\wedge_t}(\tau)=T\widehat \eta(T\tau).$$

Hence
\begin{align}
\left\|\eta\left(\tfrac{\cdot_t}T \right)F(\cdot_x,\cdot_t) \right\|^2_{X^{s,b_2}}=C\int_{\mathbb R_\xi}\langle\xi\rangle^{2s}\int_{-\infty}^{+\infty}\langle\tau+\xi^5\rangle^{2b_2}\left|\left[(T\widehat\eta(T\cdot)\ast_\tau \widehat F(\xi,\cdot)\right](\tau)\right|^2 d\tau d\xi.\label{le3_eq3}
\end{align}

Let us estimate the previous integral in the $\tau$ variable:
\begin{align}
&\notag \int_{-\infty}^{+\infty}\langle\tau+\xi^5\rangle^{2b_2}\left|\left[(T\widehat\eta(T\cdot)\ast_\tau \widehat F(\xi,\cdot)\right](\tau)\right|^2 d\tau\leq C\left\| T\widehat \eta(T\cdot)\ast_\tau\widehat F(\xi,\cdot)\right\|^2_{L^2_\tau}\\
 &+C\int_{-\infty}^{+\infty}|\tau|^{2b_2}\left|\left[(T\widehat \eta(T\cdot))\ast_\tau \widehat F(\xi,\cdot)\right](\tau-\xi^5)\right|^2d\tau\equiv I(\xi)+II(\xi).\label{le3_eq4}
\end{align}

On the one hand
\begin{align}
I(\xi)\leq C\|T\widehat \eta(T\cdot_\tau)\|^2_{L^1_\tau}\|\widehat F(\xi,\cdot_\tau)\|^2_{L^2_\tau}\leq C\|\widehat F(\xi,\cdot_\tau)\|^2_{L^2_\tau},\label{le3_eq5}
\end{align}

with $C$ independent of $T$.\\

On the other hand, to estimate $II(\xi)$, we use Leibniz formula for fractional derivatives (see \cite{KPV1993b}, Theorem A.12):
\begin{align}
\|D^\alpha(fg)-fD^\alpha g\|_{L^p(\mathbb R)}\leq C\|g\|_{L^\infty(\mathbb R)}\|D^\alpha f\|_{L^p(\mathbb R)}\quad \alpha\in(0,1),\quad 1<p<\infty.\label{le3_eq6}
\end{align}

Let us observe that 
\begin{align}
II(\xi)=C\left\|D^{b_2}_t\left(\mathcal F^{-1}_t([T\widehat\eta(T\cdot)\ast_\tau\widehat F(\xi,\cdot)](\tau-\xi^5))\right)\right\|^2_{L^2_t}.\label{le3_eq7}
\end{align}

But
\begin{align*}
\mathcal F^{-1}_t[(T\widehat \eta(T\cdot)\ast_\tau \widehat F(\xi,\cdot))(\tau-\xi^5)](t)&=C\int_{-\infty}^{+\infty}e^{it\tau}e^{it\xi^5}[T\widehat \eta(T\cdot)\ast_\tau \widehat F(\xi,\cdot)](\tau)d\tau\\
&=Ce^{it\xi^5}\mathcal F^{-1}_t[T\widehat \eta(T\cdot)](t)[F(\cdot_x,t)]^{\wedge_x}(\xi)\\
&=Ce^{it\xi^5}\eta\left(\tfrac tT\right) [F(\cdot_x,t)]^{\wedge_x}(\xi).
\end{align*}

Then, from \eqref{le3_eq7}, using \eqref{le3_eq6}, we have that
\begin{align}
\notag II(\xi)=&C\left\|D_t^{b_2}\left(e^{i\cdot_t\xi^5}[F(\cdot_x,\cdot_t)]^{\wedge_x}(\xi)\eta\left(\tfrac {\cdot_t}T\right)\right)\right\|^2_{L^2_t}\\
\notag \leq &C\left\|D_t^{b_2}\left(e^{i\cdot_t\xi^5}[F(\cdot_x,\cdot_t)]^{\wedge_x}(\xi)\eta\left(\tfrac {\cdot_t} T\right)\right)-e^{i\cdot_t\xi^5}[F(\cdot_x,\cdot_t)]^{\wedge_x}(\xi)D_t^{b_2}\left(\eta\left(\tfrac {\cdot_t}T\right)\right)\right\|^2_{L^2_t}\\
\notag&+C \left\|e^{i\cdot_t\xi^5}[F(\cdot_x,\cdot_t)]^{\wedge_x}(\xi)D_t^{b_2}\left(\eta\left(\tfrac {\cdot_t}T \right)\right) \right\|^2_{L^2_t}\\
\notag \leq & C\left\|\eta\left(\tfrac {\cdot_t}T \right) \right\|^2_{L^\infty_t}\left\|D_t^{b_2}\left(e^{i\cdot_t\xi^5}[F(\cdot_x,\cdot_t)]^{\wedge_x}(\xi) \right) \right\|^2_{L^2_t}\\
\notag &+C \left\|e^{i\cdot_t\xi^5}[F(\cdot_x,\cdot_t)]^{\wedge_x}(\xi)D_t^{b_2}\left(\eta\left(\tfrac {\cdot_t}T \right)\right) \right\|^2_{L^2_t}\\
\leq & C\left\|D_t^{b_2}\left(e^{i\cdot_t\xi^5}[F(\cdot_x,\cdot_t)]^{\wedge_x}(\xi) \right) \right\|^2_{L^2_t}+C \left\|e^{i\cdot_t\xi^5}[F(\cdot_x,\cdot_t)]^{\wedge_x}(\xi)D_t^{b_2}\left(\eta\left(\tfrac {\cdot_t}T \right)\right) \right\|^2_{L^2_t}.\label{le3_eq8}
\end{align}

From \eqref{le3_eq5} it is clear that
\begin{align}
\int_{\mathbb R_\xi}\langle\xi\rangle^{2s}I(\xi)d\xi\leq C\int_{\mathbb R_\xi}\langle\xi\rangle^{2s}\|\widehat F(\xi,\cdot_\tau)\|^2_{L^2_\tau}d\xi=C\|F\|^2_{X^{s,0}}.\label{le3_eq9}
\end{align}

Besides, from \eqref{le3_eq8}, using Plancherel's identity, we have that
\begin{align*}
\int_{\mathbb R_\xi}\langle\xi\rangle^{2s}II(\xi)d\xi\leq &C\int_{\mathbb R_\xi}\langle\xi\rangle^{2s}\int_{-\infty}^{+\infty}|\tau|^{2b_2}\left|\left(e^{it\xi^5}[F(\cdot_x,t)]^{\wedge_x}(\xi) \right)^{\wedge_t}(\tau)\right|^2 d\tau d\xi\\
&+C\int_{\mathbb R_\xi}\langle\xi\rangle^{2s}\int_{-\infty}^{+\infty}\left| e^{it\xi^5}[F(\cdot_x,t)]^{\wedge_x}(\xi)D_t^{b_2}\left(\eta\left( \tfrac tT\right)\right)\right|^2 dt d\xi\\
\leq &C\int_{\mathbb R_\xi}\langle\xi\rangle^{2s}\int_{-\infty}^{+\infty}|\tau+\xi^5|^{2b_2}|\widehat F(\xi,\tau)|^2d\tau d\xi\\
&+C\int_{\mathbb R_\xi}\langle\xi\rangle^{2s}\int_{-\infty}^{+\infty}\left| e^{it\xi^5}[F(\cdot_x,t)]^{\wedge_x}(\xi)D_t^{b_2}\left(\eta\left(\tfrac tT \right) \right)\right|^2dt d\xi.
\end{align*}

Using Hölder's inequality in the integral in the $t$ variable with conjugate exponents $p$ and $p'$, we have that
\begin{align*}
\int_{\mathbb R_\xi}\langle\xi\rangle^{2s}II(\xi)d\xi\leq C\|F\|^2_{X^{s,b_2}}+C\int_{\mathbb R_\xi}&\langle\xi\rangle^{2s}\left\{\int_{-\infty}^{+\infty}\left| e^{it\xi^5}[F(\cdot_x,t)]^{\wedge_x}(\xi)\right|^{2p}dt \right\}^{\frac1p}\\
&\left\{ \int_{-\infty}^{+\infty}|D_t^{b_2}\left(\eta\left( \tfrac tT\right) \right)|^{2p'}dt\right\}^{\frac1{p'}}d\xi.
\end{align*}

Let us choose $p$ such that $\frac12-\frac1{2p}=b_2$; i.e., $\frac1{2p}=\frac{1-2b_2}2$, which is equivalent to $2p=\frac2{1-2b_2}$, which implies $H^{b_2}(\mathbb R_t)\hookrightarrow L^{2p}(\mathbb R_t)$. Then

\begin{align*}
\int_{\mathbb R_\xi}\langle\xi\rangle^{2s}II(\xi)d\xi\leq & C\|F\|^2_{X^{s,b_2}}+C\int_{\mathbb R_\xi}\langle\xi\rangle^{2s}\left\| e^{it\xi^5}[F(\cdot_x,t)]^{\wedge_x}(\xi)\right\|^2_{H^{b_2}(\mathbb R_t)}d\xi\\
& \left\{ \int_{-\infty}^{+\infty}|D_t^{b_2}\left(\eta\left(\tfrac tT \right) \right)|^{2p'}dt\right\}^{\frac1{p'}}\\
\leq & C\|F\|^2_{X^{s,b_2}}+C\|F\|^2_{X^{s,b_2}}\left\{ \int_{-\infty}^{+\infty}|D_t^{b_2}\left(\eta\left(\tfrac tT \right) \right)|^{2p'}dt\right\}^{\frac1{p'}}.
\end{align*}

Taking into account that the inverse Fourier transform is a bounded operator from $L^{\frac{2p'}{2p'-1}}$ to $L^{2p'}$, it follows that
\begin{align}
\notag\int_{\mathbb R_\xi}\langle\xi\rangle^{2s}II(\xi)d\xi\leq & C\|F\|^2_{X^{s,b_2}}+C\|F\|^2_{X^{s,b_2}}\left\{ \int_{-\infty}^{+\infty}\left| |\tau|^{b_2}T\widehat \eta(T\tau) \right|^{\frac{2p'}{2p'-1}}d\tau\right\}^{\frac{2p'-1}{p'}}\\
\leq & C\|F\|^2_{X^{s,b_2}}.\label{le3_eq10}
\end{align}

From \eqref{le3_eq3}, \eqref{le3_eq4}, \eqref{le3_eq9}, and \eqref{le3_eq10} we conclude that
$$\left\| \eta\left(\tfrac {\cdot_t}T \right)F(\cdot_x,\cdot_t)\right\|_{X^{s,b_2}}\leq C\|F\|_{X^{s,b_2}}\quad T\in(0,1],$$

with $C$ independent of $T$.\\

In consequence of (i), and (ii), we conclude that
$$\left\| \eta\left(\tfrac{\cdot_t}T\right) F(\cdot_x,\cdot_t)\right\|_{X^{s,0}}\leq CT^{b_2}\|F\|_{X^{s,b_2}}.$$

\end{proof}

The following calculus inequalities will be used the proof of the estimative of the bilinear form (Lemma \ref{ble1}):
\begin{lemma}\label{CI}
\begin{enumerate}
\item[(i)] For $\beta\geq\gamma\geq0$, and $\quad \beta+\gamma>1$ it follows that
\begin{align}
\int_{\mathbb R}\frac1{\langle x-a_1\rangle^\beta\langle x-a_2\rangle^\gamma}dx\leq C\frac{\phi_\beta(a_1-a_2)}{\langle a_1-a_2\rangle^\gamma},\label{CI_eq1}
\end{align}

where
$$\phi_\beta(a)\sim \left\{ \begin{array}{lr}
1&\text{for }\beta>1,\\
\log(1+\langle a\rangle)&\text{for }\beta=1,\\
\langle a\rangle^{1-\beta}&\text{for }\beta<1.
\end{array}\right.
$$
(For a proof of \eqref{CI_eq1}, see \cite{ET2013}).

\item[(ii)] For $\rho\in(\tfrac12,1)$,

\begin{align}
\int_{\mathbb R}\frac1{\langle x\rangle^\rho\sqrt{|x-a|}}\leq C\frac1{\langle a\rangle^{\rho-\frac12}}.\label{CI_eq2}
\end{align}

(For a proof of \eqref{CI_eq2} see \cite{ET2016}).

\end{enumerate}
\end{lemma}

\begin{lemma}\label{ble1} (Estimative of the bilinear form). Let $s\geq 0$,  $\frac 25\leq b<\frac12$, and $0\leq a\leq 10b-4$. There exists $C>0$ such that
\begin{align}
 \|\partial_x(vw)\|_{X^{s+a,-b}}\leq C\|v\|_{X^{s,b}}\|w\|_{X^{s,b}},\label{ble1_eq1}.
\end{align}
\end{lemma}

\begin{proof} We observe that
$$[\partial_x(vw)]^\wedge(\xi,\tau)=Ci\xi\int_{\mathbb R_2} \widehat v(\xi_1,\tau_1)\widehat w(\xi-\xi_1,\tau-\tau_1)d\xi_1 d\tau_1.$$

Hence
\begin{align}
\notag &\|\partial_x(vw)\|^2_{X^{s+a,-b}}\\
\notag&=C\int_{\mathbb R^2_{\xi\tau}}\langle\xi\rangle^{2(s+a)}\langle\tau+\xi^5\rangle^{-2b}\xi^2\left|\int_{\mathbb R^2_{\xi_1\tau_1}}\widehat v(\xi_1,\tau_1)\widehat w(\xi-\xi_1,\tau-\tau_1)d\xi_1 d\tau_1 \right|^2d\xi d\tau\\
&=C\int_{\mathbb R^2_{\xi\tau}}\left| \int_{\mathbb R^2_{\xi_1\tau_1}} \langle\xi\rangle^{(s+a)}\langle\tau+\xi^5\rangle^{-b}|\xi|\widehat v(\xi_1,\tau_1)\widehat w(\xi-\xi_1,\tau-\tau_1)d\xi_1 d\tau_1 \right|^2d\xi d\tau.\label{ble1_eq2}
\end{align}

Let $h\in L^2(\mathbb R^2_{\xi\tau})$ an arbitrary function. If we manage to prove that
\begin{align}
\notag&\left|\int_{\mathbb R^2_{\xi\tau}} \left[\int_{\mathbb R^2_{\xi_1\tau_1}}\langle\xi\rangle^{(s+a)}\langle\tau+\xi^5\rangle^{-b}|\xi|\widehat v(\xi_1,\tau_1)\widehat w(\xi-\xi_1,\tau-\tau_1)d\xi_1 d\tau_1 \right]h(\xi,\tau)d\xi d\tau\right|\\
&\leq C\|v\|_{X^{s,b}}\|w\|_{X^{s,b}}\|h\|_{L^2(\mathbb R^2)},\label{ble1_eq3}
\end{align}

then we would have, by a duality argument, that
\begin{align}
\|\partial_x(vw)\|_{X^{s+a,-b}}\leq C\|v\|_{X^{s,b}}\|w\|_{X^{s,b}}.\label{ble1_eq4}
\end{align}

Taking into account that, there exists $C>0$, such that for $s\geq 0$
$$\frac{\langle\xi\rangle^s}{\langle\xi_1\rangle^s\langle\xi-\xi_1\rangle^s}\leq C,$$
then, to establish \eqref{ble1_eq3}, it is enough to prove that
\begin{align}
\notag \Bigg |\int_{\mathbb R^2_{\xi\tau}} \Bigg[\int_{\mathbb R^2_{\xi_1\tau_1}}& \langle\xi\rangle^{a}\langle\xi_1\rangle^s\langle\tau+\xi^5\rangle^{-b}|\xi|\langle\xi-\xi_1\rangle^s|\widehat v(\xi_1,\tau_1)||\widehat w(\xi-\xi_1,\tau-\tau_1)|d\xi_1 d\tau_1 \Bigg]\\
&h(\xi,\tau)d\xi d\tau\Bigg|\leq C\|v\|_{X^{s,b}}\|w\|_{X^{s,b}}\|h\|_{L^{2}(\mathbb R^2)}.\label{ble1_eq5}
\end{align}

Since
\begin{align*}
\Bigg|\int_{\mathbb R^2_{\xi\tau}} \Bigg[\int_{\mathbb R^2_{\xi_1\tau_1}} &\langle\xi\rangle^a\langle\xi_1\rangle^s\langle\tau+\xi^5\rangle^{-b}|\xi|\langle\xi-\xi_1\rangle^s|\widehat v(\xi_1,\tau_1)||\widehat w(\xi-\xi_1,\tau-\tau_1)|d\xi_1 d\tau_1 \Bigg]\\
& h(\xi,\tau)d\xi d\tau\Bigg|\\
\leq\hspace{2mm} \int_{\mathbb R^2_{\xi\tau}} \int_{\mathbb R^2_{\xi_1\tau_1}} &\frac{\langle\xi\rangle^a|\xi||h(\xi,\tau)| \langle\xi_1\rangle^s\langle\tau_1+\xi_1^5\rangle^{b}|\widehat v(\xi_1,\tau_1)|\langle\xi-\xi_1\rangle^s\langle\tau-\tau_1+(\xi-\xi_1)^5\rangle^b}{\langle\tau+\xi^5\rangle^b\langle\tau_1+\xi_1^5\rangle^b\langle\tau-\tau_1+(\xi-\xi_1)^5\rangle^b}\\
&|\widehat w(\xi-\xi_1,\tau-\tau_1)|d\xi_1d\tau_1 d\xi d\tau\\
\leq\int_{\mathbb R^2_{\xi\tau}}\Bigg[ \int_{\mathbb R^2_{\xi_1\tau_1}}&\frac{\langle\xi\rangle^{2a}|\xi|^2|h(\xi,\tau)|^2}{\langle\tau+\xi^5\rangle^{2b}\langle\tau_1+\xi_1^5\rangle^{2b}\langle\tau-\tau_1+(\xi-\xi_1)^5\rangle^{2b}}d\xi_1 d\tau_1\Bigg]^{\frac12}\\
\Bigg[\int_{\mathbb R^2_{\xi_1\tau_1}}&\langle\xi_1\rangle^{2s}\langle\tau_1+\xi_1^5\rangle^{2b}|\widehat v(\xi_1,\tau_1)|^2\langle\xi-\xi_1\rangle^{2s}\langle\tau-\tau_1+(\xi-\xi_1)^5\rangle^{2b}\\
&|\widehat w(\xi-\xi_1,\tau-\tau_1)|^2d\xi_1 d\tau_1 \Bigg]^{\frac12}d\xi d\tau\\
\leq\Bigg[\int_{\mathbb R^2_{\xi\tau}}\int_{\mathbb R^2_{\xi_1\tau_1}}&\frac{\langle\xi\rangle^{2a}|\xi|^2|h(\xi,\tau)|^2}{\langle\tau+\xi^5\rangle^{2b}\langle\tau_1+\xi_1^5\rangle^{2b}\langle\tau-\tau_1+(\xi-\xi_1)^5\rangle^{2b}}d\xi_1 d\tau_1 d\xi d\tau\Bigg]^{\frac12}\\
\Bigg[\int_{\mathbb R^2_{\xi\tau}}\int_{\mathbb R^2_{\xi_1\tau_1}}&\langle\xi_1\rangle^{2s}\langle\tau_1+\xi_1^5\rangle^{2b}|\widehat v(\xi_1,\tau_1)|^2\langle\xi-\xi_1\rangle^{2s}\langle\tau-\tau_1+(\xi-\xi_1)^5\rangle^{2b}\\
&|\widehat w(\xi-\xi_1,\tau-\tau_1)|^2d\xi_1 d\tau_1d\xi d\tau \Bigg]^{\frac12}
\end{align*}
\begin{align*}
=\Bigg[ \int_{\mathbb R^2_{\xi\tau}}&|h(\xi,\tau)|^2\left(\int_{\mathbb R^2_{\xi_1\tau_1}}\frac{\langle\xi\rangle^{2a}|\xi|^2}{\langle\tau+\xi^5\rangle^{2b}\langle\tau_1+\xi_1^5\rangle^{2b}\langle\tau-\tau_1+(\xi-\xi_1)^5\rangle^{2b}}d\xi_1 d\tau_1 \right)d\xi d\tau\Bigg]^{\frac12}\\
&\|v\|_{X^{s,b}}\|w\|_{X^{s,b}},
\end{align*}

then, to prove \eqref{ble1_eq5}, it is enough to prove that
\begin{align}
\sup_{(\xi,\tau)\in\mathbb R^2}\int_{\mathbb R^2_{\xi_1\tau_1}}\frac{\langle\xi\rangle^{2a}|\xi|^2}{\langle\tau+\xi^5\rangle^{2b}\langle\tau_1+\xi_1^{5}\rangle^{2b}\langle\tau-\tau_1+(\xi-\xi_1)^5\rangle^{2b}}d\xi_1d\tau_1\leq C.\label{ble1_eq6}
\end{align}

Let us observe that
\begin{align*}
&\int_{\mathbb R_{\tau_1}}\frac{\langle\xi\rangle^{2a}|\xi|^2}{\langle\tau+\xi^5\rangle^{2b}\langle\tau_1+\xi_1^{5}\rangle^{2b}\langle\tau-\tau_1+(\xi-\xi_1)^5\rangle^{2b}}d\tau_1\\
&=\frac{\langle\xi\rangle^{2a}|\xi|^2}{\langle\tau+\xi^5\rangle^{2b}}\int_{\mathbb R_{\tau_1}}\frac1{\langle\tau_1-(-\xi_1^5)\rangle^{2b}\langle\tau_1-(\tau+(\xi-\xi_1)^5)\rangle^{2b}}d\tau_1.
\end{align*}

Using inequality \eqref{CI_eq1}, with $\beta=\gamma=2b<1$, and $\beta+\gamma=4b>1$ $(b>\frac14)$, we conclude that
\begin{align}
\notag &\int_{\mathbb R_{\xi_1}}\int_{\mathbb R_{\tau_1}}\frac{\langle\xi\rangle^{2a}|\xi|^2}{\langle\tau+\xi^5\rangle^{2b}\langle\tau_1+\xi_1^5\rangle^{2b}\langle\tau-\tau_1+(\xi-\xi_1)^5\rangle^{2b}}d\tau_1d\xi_1\\
&\leq C \int_{\mathbb R_{\xi_1}}\frac{\langle\xi\rangle^{2a}|\xi|^2}{\langle\tau+\xi^5\rangle^{2b}\langle\xi_1^5+\tau+(\xi-\xi_1)^5\rangle^{4b-1}}d\xi_1.\label{ble1_eq7}
\end{align}

Let us make the following change of variables in the integral of the right hand side of \eqref{ble1_eq7}:
\begin{align}
\notag \mu&:=\xi_1^5+\tau+(\xi-\xi_1)^5=\tau+\xi^5-\frac52\xi\xi_1(\xi-\xi_1)[\xi^2+\xi_1^2+(\xi-\xi_1)^2]\\
&=\tau+\frac1{16}\xi^5+\frac52\xi^3(\xi_1-\tfrac\xi2)^2+5\xi(\xi_1-\tfrac\xi2)^4,\label{ble1_eq8}
\end{align}

\begin{align}
d\mu=-\frac52\xi(\xi-2\xi_1)[\xi^2+(\xi-2\xi_1)^2]d\xi_1.\label{ble1_eq9}
\end{align}

Therefore
\begin{align}
\notag\int_{\mathbb R_{\xi_1}}\int_{\mathbb R_{\tau_1}}&\frac{\langle\xi\rangle^{2a}|\xi|^2}{\langle\tau+\xi^5\rangle^{2b}\langle\tau_1+\xi_1^5\rangle^{2b}\langle\tau-\tau_1+(\xi-\xi_1)^5\rangle^{2b}}d\tau_1d\xi_1\\
&\leq C\frac{\langle\xi\rangle^{2a}|\xi|^2}{\langle\tau+\xi^5\rangle^{2b}}\int_{A\subset\mathbb R_{\mu}}\frac{1}{(1+|\mu|)^{4b-1}|\xi||\xi-2\xi_1|[\xi^2+(\xi-2\xi_1)^2]}d\mu,\label{ble1_eq10}
\end{align}

where $A=[\tau+\tfrac1{16}\xi^5,+\infty)$ when $\xi>0$.\\

Using the last equality of \eqref{ble1_eq8}, it can be seen that
$$\left(\xi_1-\tfrac\xi2 \right)^2=-\frac{\xi^2}4+\sqrt{\frac{\xi^5-4\tau+4\mu}{20\xi}},$$
then
\begin{align*}
(\xi-2\xi_1)^2&=-\xi^2+\sqrt{\frac4{5\xi}(\xi^5-4\tau+4\mu)},\\
\xi_1&=\frac\xi2\pm\sqrt{\frac12\sqrt{\frac1{5\xi}(\xi^5-4\tau+4\mu)}-\frac14\xi^2},\\
2\xi_1-\xi&=\pm\sqrt{-\xi^2+\sqrt{\frac4{5\xi}(\xi^5-4\tau+4\mu)}}.
\end{align*}

Therefore, from \eqref{ble1_eq10}, we have that
\begin{align}
\notag &\int_{\mathbb R^2_{\xi_1\tau_1}}\frac{\langle\xi\rangle^{2a}|\xi|^2}{\langle\tau+\xi^5\rangle^{2b}\langle\tau_1+\xi_1^5\rangle^{2b}\langle\tau-\tau_1+(\xi-\xi_1)^5\rangle^{2b}}d\tau_1 d\xi_1\\
\notag&\leq C\frac{\langle\xi\rangle^{2a}|\xi|}{\langle\tau+\xi^5\rangle^{2b}}\int_{A\subset \mathbb R_\mu}\frac1{(1+|\mu|)^{4b-1}\sqrt{-\xi^2+\sqrt{\frac4{5\xi}(\xi^5-4\tau+4\mu)}}\sqrt{\frac4{5\xi}(\xi^5-4\tau+4\mu)}}d\mu\\
&\leq C\frac{\langle\xi\rangle^{2a}|\xi|}{\langle\tau+\xi^5\rangle^{2b}}\int_{A\subset \mathbb R_\mu}\frac{\sqrt{\xi^2+\sqrt{\frac4{5\xi}(\xi^5-4\tau+4\mu)}}}{(1+|\mu|)^{4b-1}\sqrt{\frac{-16\tau+16\mu-\xi^5}{5\xi}}\sqrt{\frac4{5\xi}(\xi^5-4\tau+4\mu)}}d\mu.\label{ble1_eq11}
\end{align}

Let us assume, without loss of generality, that $\xi>0$, and let us observe that the function
$$\mu:=\varphi(\xi_1)=\tau+\xi_1^5+(\xi-\xi_1)^5$$
has its absolute minimum in $\xi_1=\frac\xi2$. Then, the integration interval $A$ in \eqref{ble1_eq11} is
$$A=[\varphi(\tfrac\xi2),+\infty)=[\tau+\tfrac1{16}\xi^5,+\infty).$$

This way
$$\mu\geq\tau+\frac1{16}\xi^5,\quad \xi^5-4\tau+4\mu\geq \xi^5-4\tau+4\tau+\frac14\xi^5=\frac54\xi^5,$$
\begin{align}
\frac4{5\xi}(\xi^5-4\tau+4\mu)&\geq\frac4{5\xi}\frac54\xi^5=\xi^4,\label{ble1_eq12}\\
\sqrt{\frac4{5\xi}(\xi^5-4\tau+4\mu)}&\geq\xi^2.\label{ble1_eq13}
\end{align}

Hence, taking into account \eqref{ble1_eq13}, \eqref{ble1_eq12}, from \eqref{ble1_eq11}, we conclude that
\begin{align}
\notag &\int_{\mathbb R^2_{\xi_1\tau_1}}\frac{\langle\xi\rangle^{2a}|\xi|^2}{\langle\tau+\xi^5\rangle^{2b}\langle\tau_1+\xi_1^5\rangle^{2b}\langle\tau-\tau_1+(\xi-\xi_1)^5\rangle^{2b}}d\tau_1 d\xi_1\\
\notag &\leq C\frac{\langle\xi\rangle^{2a}|\xi|}{\langle\tau+\xi^5\rangle^{2b}}\int_{\mathbb R_\mu}\frac1{\langle \mu\rangle^{4b-1}\sqrt{\frac{-16\tau+16\mu-\xi^5}{5\xi}}|\xi|}d\mu\\
&\leq C\frac{\langle\xi\rangle^{2a}|\xi|^{\frac12}}{\langle\tau+\xi^5\rangle^{2b}}\int_{\mathbb R_\mu}\frac1{\langle \mu\rangle^{4b-1}\sqrt{|\mu-\tau-\frac1{16}\xi^5|}}d\mu.\label{ble1_eq14}
\end{align}

Using the calculus inequality \eqref{CI_eq2} (Lemma \ref{CI}) with $\rho:=4b-1$ ($4b-1\in(\tfrac12,1)\iff \frac38<b<\frac12$) in \eqref{ble1_eq14}, it follows that
\begin{align}
\int_{\mathbb R^2_{\xi_1\tau_1}}\frac{\langle\xi\rangle^{2a}|\xi|^2}{\langle\tau+\xi^5\rangle^{2b}\langle\tau_1+\xi_1^5\rangle^{2b}\langle\tau-\tau_1+(\xi-\xi_1)^5\rangle^{2b}}d\tau_1 d\xi_1\leq C\frac{\langle\xi\rangle^{2a}|\xi|^{\frac12}}{\langle\tau+\xi^5\rangle^{2b}\langle\tau+\tfrac1{16}\xi^5\rangle^{4b-\frac32}}.\label{ble1_eq15}
\end{align}

It is clear, from \eqref{ble1_eq15}, for $|\xi|\leq 1$, and $\frac38< b<\frac12$, that
\begin{align}
\int_{\mathbb R^2_{\xi_1\tau_1}}\frac{\langle\xi\rangle^{2a}|\xi|^2}{\langle\tau+\xi^5\rangle^{2b}\langle\tau_1+\xi_1^5\rangle^{2b}\langle\tau-\tau_1+(\xi-\xi_1)^5\rangle^{2b}}d\tau_1 d\xi_1\leq C.\label{ble1_eq16}
\end{align}
Let us assume then that $|\xi|>1$. We observe that
$$|\xi^5|\leq \frac{16}{15}|\tau+\xi^5|+\frac{16}{15}|\tau+\tfrac1{16}\xi^5|,$$
then
\begin{align*}
\frac12|\xi^5|\leq \frac{16}{15}|\tau+\xi^5|\quad \text{or}\quad \frac12|\xi^5|\leq \frac{16}{15}|\tau+\tfrac1{16}\xi^5|.
\end{align*}

For $\frac12|\xi^5|\leq \frac{16}{15}|\tau+\xi^5|$, we have that
\begin{align}
C\frac{\langle\xi\rangle^{2a}|\xi|^{\frac12}}{\langle\tau+\xi^5\rangle^{2b}\langle\tau+\tfrac1{16}\xi^5\rangle^{4b-\frac32}}\leq C\frac{|\xi|^{2a+\frac12}}{\langle\tau+\xi^5\rangle^{2b}}\leq C\frac{|\xi|^{2a+\frac12}}{\langle\tfrac{15}{32}|\xi^5|\rangle^{2b}}\leq C|\xi|^{2a+\frac12-10b}\leq C,\label{ble1_eq17}
\end{align}

since $2a+\frac12-10b\leq20b-8+\frac12-10b=10b-\frac{15}2<0$.\\

For $\frac12|\xi^5|\leq\frac{16}{15}|\tau+\tfrac1{16}\xi^5|$, we have that
\begin{align}
\notag C\frac{\langle\xi\rangle^{2a}|\xi|^{\frac12}}{\langle\tau+\xi^5\rangle^{2b}\langle\tau+\tfrac1{16}\xi^5\rangle^{4b-\frac32}}&\leq C\frac{|\xi|^{2a+\frac12}}{\langle\tau+\tfrac1{16}\xi^5\rangle^{4b-\frac32}}\leq C\frac{|\xi|^{2a+\frac12}}{\langle\tfrac{15}{32}|\xi^5|\rangle^{4b-\frac32}}\leq C|\xi|^{2a+\frac12-20b+\frac{15}2}\\
&\leq C,\label{ble1_eq18}
\end{align}
since $2a+\frac12-20b+\frac{15}2\leq20b-8+\frac12-20b+\frac{15}2=0$.\\

In this way, from \eqref{ble1_eq15} to \eqref{ble1_eq18}, it follows that, for $\frac25<b<\frac12$ and $0\leq a\leq 10b-4$, inequality \eqref{ble1_eq6} holds, which proves \eqref{ble1_eq4}.
\end{proof}

\section{Existence of a local solution in time of IVP \eqref{maineq}. Proof of Theorem \ref{localsolintime}}

Let us demonstrate that there exists $T\in (0,\frac12]$ such that the operator $\Gamma_T$, defined for $u\in X^{s,b}$ by
\begin{align}
(\Gamma_T u)(t) := \eta(t) W(t) u_0 + \eta(t) \int_0^T W(t-t') F_T(u(t')) dt', \label{gamma}
\end{align}
has a fixed point in $X^{s,b}$.\\

Recall that by inequality \eqref{le1_eq1} in Lemma \ref{le1}, we have that
\begin{align}
\|\eta(\cdot_t)[W (\cdot_t)u_0](\cdot_x)\|_{X^{s,b}(\mathbb R^2)}\leq C\|u_0\|_{H^s(\mathbb R_x)},\label{lst_eq1}
\end{align}

and by using Lemma \ref{le2} with $T=1$, $b':=-b^*$, where $b^*$ is such that $b<b^*<\frac12$,
\begin{align}
\left\|\eta (\cdot_t)\int_0^{\cdot_t} W (\cdot_t{\scriptstyle -}t')F_T(u(t')) dt'\right\|_{X^{s,b}}\leq C \|F_T(u)\|_{X^{s,-b^*}}.\label{lst_eq2}
\end{align}

Since $-\frac12<-b^*<-b<\frac12$, by Lemma \ref{le3}, with $2T\leq 1$, and Lemma \ref{ble1}, we have that
\begin{align}
\|F_T(u)\|_{X^{s,-b^*}} & \leq CT^{-b+b^*} \|-\frac12 \partial_x(u^2)\|_{X^{s,-b}}\leq C T^{-b+b^*} \|u\|^2_{X^{s,b}}.\label{lst_eq3}
\end{align}

Therefore, combining \eqref{lst_eq2} and \eqref{lst_eq3}, we can affirm that
\begin{align}
\left\|\eta (\cdot_t)\int_0^{\cdot_t} W (\cdot_t{\scriptstyle -}t')F_T(u(t')) dt'\right\|_{X^{s,b}}\leq C T^{-b+b^*} \|u\|^2_{X^{s,b}}.\label{lst_eq4}
\end{align}

From \eqref{lst_eq1} and \eqref{lst_eq4} we obtain
\begin{align}
\|\Gamma_T u \|_{X^{s,b}} \leq C (\|u_0\|_{H^s(\mathbb R_x)} + T^{-b+b^*} \|u\|^2_{X^{s,b}}). \label{lst_eq5}
\end{align}

Let us define the closed ball $B_{X^{s,b}}(0,R)\subset X^{s,b}$ centered at 0, with radius $R:=2C\|u_0\|_{H^s(\mathbb R_x)}>0$. Hence, from \eqref{lst_eq5}, it follows that if $u\in B_{X^{s,b}}(0,R)$,
$$\|\Gamma_T u \|_{X^{s,b}} \leq \frac R2 + C T^{-b+b^*} R^2.$$

Choosing $T>0$ such that
\begin{align}
CT^{-b+b^*} R^2 \leq \frac R2, \label{lst_eq6}
\end{align}

we have that $\Gamma_T$ maps $B_{X^{s,b}}(0,R)$ into itself.\\

It remains to prove that $\Gamma_T: B_{X^{s,b}}(0,R) \to B_{X^{s,b}}(0,R)$ is a contraction. Let us consider $v,w\in B_{X^{s,b}}(0,R)$. Then, by using Lemmas \ref{le2}, \ref{le3}, and \ref{ble1}, we have that
\begin{align}
\notag \|\Gamma_T v - \Gamma_T w \|_{X^{s,b}} & = \left\|\eta(\cdot_t)\int_0^{\cdot_t} W (\cdot_t{\scriptstyle -}t') [F_T(v(t')) - F_T(w(t'))] dt'\right\|_{X^{s,b}}\\
\notag& \leq  CT^{-b+b^*} \|\partial_x(v^2 - w^2)\|_{X^{s,-b}}\\
& \leq C T^{-b+b^*} \|v+w\|_{X^{s,b}} \|v-w\|_{X^{s,b}}\leq 2 CR T^{-b+b^*} \|v-w\|_{X^{s,b}}. \label{lst_eq7}
\end{align}

If we choose $T>0$ such that $2CRT^{-b+b^*} <1$, it is clear that $T$ also satisfies \eqref{lst_eq6}. In consequence, from \eqref{lst_eq7} we conclude that $\Gamma_T:B_{X^{s,b}}(0,R) \to B_{X^{s,b}}(0,R)$ is a contraction.\\

Let us observe that we can also assume that $T\in (0,\frac12]$. Therefore, there exist $T\in(0,\frac12]$ and a unique $u\in B_{X^{s,b}}(0,R)$ such that $\Gamma_T u=u$.\\

Now we prove that $u\in C(\mathbb R_t; H^s(\mathbb R_x))$. Let us note that the first term on the right hand side of \eqref{gamma} is continuous from $\mathbb R_t$, with values in $H^s(\mathbb R_x)$. To demonstrate that the second term on the right hand side of \eqref{gamma} is continuous form $\mathbb R_t$, with values in $H^s(\mathbb R_x)$, it is enough to see that this term belongs to $X^{s,\widetilde b}$ for some $\widetilde b>\frac12$.\\

Let us take $\widetilde b>\frac12$ such that $b^*+\widetilde b\leq 1$. Then
$$-\frac12<-b^*<0\leq \widetilde b\leq (-b^*)+1.$$
By applying Lemma \ref{le2} with $T=1$, and Lemmas \ref{le3} and \ref{ble1}, we conclude that
\begin{align*}
\left\|\eta(\cdot_t)\int_0^{\cdot_t} W (\cdot_t{\scriptstyle -}t')F_T(u(t')) dt'\right\|_{X^{s,\widetilde b}}&\leq C \|\eta \left(\tfrac{\cdot_t}{2T} \right) \partial_x u^2\|_{X^{s,-b^*}} \leq C T^{-b+b^*} \|\partial_x u^2\|_{X^{s,-b}}\\
& \leq C T^{-b+b^*}\|u\|^2_{X^{s,b}} <\infty.
\end{align*}

In this way, we have proved that $u=\Gamma_T u \in C(\mathbb R_t; H^s(\mathbb R_x))$. Theorem \ref{localsolintime} is proved.\qed

\section{Regularity gain of the nonlinear part of the solution of IVP \eqref{maineq}. Proof of Theorem \ref{more_reg}}

Let $b^*$ and $\widetilde b$ be such that $b<b^*<\frac12 < \widetilde b$, and $b^*+\widetilde b\leq 1$. By using the immersion $X^{s+a,\widetilde b} \hookrightarrow C_b(\mathbb R_t; H^{s+a}(\mathbb R_x))$ and taking into account that for $0\leq t'\leq t\leq T\leq\frac12$, $\eta(t)=1$ and $-\frac12 \partial_x (u(t'))^2=F_T(u(t'))$, we have that
\begin{align}
\notag&\left \|\int_0^{\cdot_t} W (\cdot_t{\scriptstyle -}t') (-\frac12 \partial_x(u(t'))^2) dt' \right\|_{C([0,T];H^{s+a}(\mathbb R_x))}\\
\notag& = \left\| \eta(\cdot_t) \int_0^{\cdot_t} W (\cdot_t{\scriptstyle -}t')F_T(u(t')) dt'  \right\|_{C([0,T];H^{s+a}(\mathbb R_x))}\\
\notag &\leq \left\| \eta(\cdot_t) \int_0^{\cdot_t} W (\cdot_t{\scriptstyle -}t')F_T(u(t')) dt'  \right\|_{C_b(\mathbb R_t;H^{s+a}(\mathbb R_x))}\\
& \leq C\left\| \eta(\cdot_t) \int_0^{\cdot_t} W (\cdot_t{\scriptstyle -}t')F_T(u(t')) dt'  \right\|_{X^{s+a,\tilde b}}.\label{PT1_2-1}
\end{align}

Applying Lemma \ref{le2} with $-\frac12<-b^*\leq 0\leq \tilde b\leq -b^*+1$, it follows that
\begin{align}
\left\| \eta(\cdot_t) \int_0^{\cdot_t} W (\cdot_t{\scriptstyle -}t')F_T(u(t')) dt'  \right\|_{X^{s+a,\tilde b}} \leq C \|F_T(u(\cdot))\|_{X^{s+a,-b^*}}.\label{PT1_2-2}
\end{align}

From Lemma \ref{le3} with $b_1:=-b^*$, $b_2:=-b$, and Lemma \ref{ble1}, we can conclude that

\begin{align}
\| F_T(u(\cdot)) \|_{X^{s+a,-b^*}} \leq CT^{-b+b^*} \|\partial_x u^2 \|_{X^{s+a,-b}} \leq C \|u\|^2_{X^{s,b}}. \label{PT1_2-3}
\end{align}

Finally, from \eqref{PT1_2-1}, \eqref{PT1_2-2} and \eqref{PT1_2-3}, we have that
\begin{align*}
\left \|\int_0^{\cdot_t} W (\cdot_t{\scriptstyle -}t') (-\frac12 \partial_x(u(t'))^2) dt' \right\|_{C([0,T];H^{s+a}(\mathbb R))}\leq C \|u\|^2_{X^{s,b}}<+\infty,
\end{align*}

and Theorem \ref{more_reg} is proved.\qed

\section{Dispersive blow-up. Proof of Theorem \ref{blowup}}

Let us define $\phi(x):=e^{-2|x|}$, $x\in\mathbb R$. It is clear that $\phi\in C^1(\mathbb R\setminus \{0\})\setminus C^1(\mathbb R)$, $e^x\phi\in L^2(\mathbb R)$, and $e^{-x}\phi\in L^2(\mathbb R)$. Let $\{\alpha_j\}_{j=1}^\infty$ be a sequence of positive real numbers such that $\sum_{j=1}^\infty \alpha_j e^{4j}<+\infty$, and let us define
\begin{align}
u_0(x):=\sum_{j=1}^\infty \alpha_j [W(-j) \phi](x). \label{blowup_eq1}
\end{align}
Then we have that $u_0\in H^1(\mathbb R)$ (in fact, $u_0\in H^{\frac32-}(\mathbb R)$). Actually,
\begin{align*}
\|W(-j)\phi\|_{H^{\frac32-}} & = \left(\int_{-\infty}^{+\infty} (1 + \xi^2)^{\frac32-} |[W(-j) \phi]^\wedge (\xi) |^2 d\xi \right)^{\frac12}\\
& =\left( \int_{-\infty}^{+\infty} (1+\xi^2)^{\frac32-} |e^{ij\xi^5} \widehat \phi (\xi) |^2 d\xi \right)^{\frac12}\\
& =\left( \int_{-\infty}^{+\infty} (1+\xi^2)^{\frac32-} |\widehat \phi(\xi)|^2 d\xi \right)^{\frac12}.
\end{align*}

Since
\begin{align*}
\widehat \phi(\xi) = \frac1{\sqrt{2\pi}} \int_{-\infty}^{+\infty} e^{-i\xi x} e^{-2|x|} dx = \frac1{\sqrt{2\pi}} \left[\frac1{2+i\xi}+ \frac1{2-i\xi} \right]=\frac4{\sqrt{2\pi} (\xi^2+4)},
\end{align*}

then
\begin{align}
\|W(-j)\phi\|_{H^{\frac32-}} = \left(\int_{-\infty}^{+\infty}(1+\xi^2)^{\frac32-} \frac{16}{2\pi (\xi^2+4)^2 } d\xi \right)^{\frac12} < +\infty, \label{blowup_eq2}
\end{align}

From \eqref{blowup_eq1} and \eqref{blowup_eq2} it can be seen that $u_0\in H^{\frac32-}$, given that
\begin{align*}
\sum_{j=1}^\infty \alpha_j \|W(-j) \phi\|_{H^{\frac32-}} &= \left( \int_{-\infty}^{+\infty} (1+\xi^2)^{\frac32-} \frac{16}{2\pi (\xi^2+4)^2} d\xi \right) \sum_{j=1}^\infty \alpha_j\\
& \leq \left( \int_{-\infty}^{+\infty} (1+\xi^2)^{\frac32-} \frac{16}{2\pi (\xi^2+4)^2} d\xi \right) \sum_{j=1}^\infty \alpha_j e^{4j}<+\infty.
\end{align*}

It can be seen that the solution $u$ of the IVP \eqref{maineq} is global in time and is such that
\begin{align}
u(t) = W(t) u_0 + \int_0^t W(t-t')(-\frac12 \partial_x (u(t'))^2) dt'\equiv W(t) u_0 + z(t). \label{blowup_eq3}
\end{align}

By Theorem \ref{more_reg}, $z(t)\in C^1(\mathbb R)$ for every $t\geq 0$. (In fact, there exists $\beta>0$ such that $z(t)\in H^{(\frac32-)+\beta}$. Therefore, by taking $\frac32-$ close enough to $\frac32$ from the left, it follows that $(\frac32-)+\beta>\frac32$, and as a result $z(t)\in C^1(\mathbb R)$).\\

Consequently, the proof of Theorem \ref{blowup} reduces to establish statements (i) and (ii) for $W(t)u_0$ instead of $u(\cdot,t)$. Let us define
\begin{align}
w(x,t):=e^x[W(t)\phi](x).\label{blowup_eq3}
\end{align}

Then
\begin{align}
\left. \begin{array}{rl}
\partial_t (e^{-x }w) + \partial_x^5(e^{-x}w)&\hspace{-2mm}=0,\\
w(x,0)&\hspace{-2mm}=e^x \phi(x),
\end{array} \right\}\label{blowup_eq4}
\end{align}

which means that $w$ satisfies the following differential equation:
$$\partial_t w + \partial_x^5 w - 5\partial_x^4 w + 10 \partial_x^3 w -10 \partial_x^2 w + 5 \partial_x w -w =0,$$

which indicates that $w$ satisfies the IVP
\begin{align}
\left. \begin{array}{rl}
\partial_t w + (\partial_x -1)^5 w&\hspace{-2mm}=0,\\
w(x,0)&\hspace{-2mm}=e^x \phi(x).
\end{array} \right\}\label{blowup_eq5a}
\end{align}

Taking the Fourier transform with respect to the variable $x$ in the differential equation in \eqref{blowup_eq5a}, we obtain the IVP
\begin{align}
\left. \begin{array}{rl}
\frac d{dt} \widehat {w(\cdot,t)} + (i\xi-1)^5 \widehat {w(\cdot,t)} &\hspace{-2mm}=0,\\
\widehat{w(\cdot,0)}&\hspace{-2mm}=[e^x \phi]^\wedge(\xi),
\end{array} \right\}\label{blowup_eq5}
\end{align}

whose solution is
\begin{align}
\widehat{w(\cdot,t)} (\xi) = e^{-i\xi^5 t} e^{5t\xi^4}e^{10it\xi^3}e^{-10t\xi^2}e^{-5it \xi} e^t (e^x \phi )^\wedge (\xi). \label{blowup_eq6}
\end{align}

Let us consider
\begin{align}
\widetilde w(x,t):= e^{-x} [W(t) \phi](x). \label{blowup_eq7}
\end{align}

In a similar way as we obtained \eqref{blowup_eq6}, it can be seen that
\begin{align}
\widehat{\widetilde w(\cdot,t)} (\xi) = e^{-i\xi^5 t} e^{-5t\xi^4}e^{10it\xi^3}e^{10t\xi^2}e^{-5it \xi} e^{-t} (e^{-x} \phi )^\wedge (\xi). \label{blowup_eq8}
\end{align}

Next, we prove that $u_0\in C^\infty(\mathbb R)$. For that, it is sufficient to prove that
$$\sum_{j=1}^\infty \alpha_j e^x W(-j) \phi \in C^{\infty}(\mathbb R),$$
which will be true if for every $m\in\mathbb N\cup \{0\}$, we have
$$\sum_{j=1}^{\infty} \alpha_j \|e^x W(-j)\phi\|_{H^m(\mathbb R)}<+\infty.$$

For $m=0$, using  Plancherel's Theorem and \eqref{blowup_eq6}, it results that
\begin{align}
\notag\sum_{j=1}^\infty \alpha_j \|e^x W(-j) \phi\|_{L^2(\mathbb R)} &=\sum_{j=1}^\infty \alpha_j \|w(\cdot_x,-j)\|_{L^2_x} =\sum_{j=1}^\infty \|\widehat{w(\cdot_x,-j)}(\cdot_\xi)\|_{L^2_\xi}\\
\notag& = \sum_{j=1}^\infty \alpha_j \|e^{-5j(\xi^4-2\xi^2)}e^{-j} (e^x\phi)^\wedge(\xi)\|_{L^2_\xi}\\
\notag& = \sum_{j=1}^\infty \alpha_j \|e^{-5j(\xi^4-2\xi^2+1)} e^{4j} (e^x \phi)^\wedge(\xi)\|_{L^2_\xi}\\
&\leq \left( \sum_{j=1}^\infty \alpha_j e^{4j} \right) \|e^x \phi\|_{L^2_x}<\infty. \label{blowup_eq9}
\end{align}

Now, we prove that, for every $m\in\mathbb N$,
\begin{align}
\sum_{j=1}^\infty \alpha_j \|\partial_x^m (e^x W(-j) \phi)\|_{L^2_x} <+\infty.\label{blowup_eq10}
\end{align}

In fact,
\begin{align*}
\sum_{j=1}^\infty \alpha_j \|\partial_x^m(e^x W(-j)\phi)\|_{L^2_x} & = \sum_{j=1}^\infty \alpha_j \|\xi^m e^{-5j(\xi^2-1)^2}e^{4j}(e^x\phi)^\wedge (\xi)\|_{L^2_\xi}\\
& \leq \sum_{j=1}^\infty \alpha_j e^{4j}\|\xi^m e^{-5j(\xi^2-1)^2}\|_{L^\infty_\xi} \|e^x\phi\|_{L^2_x}\\
& \leq \sum_{j=1}^\infty \alpha_j e^{4j} \frac{\left(1+\sqrt{1+\frac m{5j}} \right)^{m/2}}{2^{m/2}} \|e^x\phi\|_{L^2_x}\\
& \leq  \frac{\left(1+\sqrt{1+\frac m{5}} \right)^{m/2}}{2^{m/2}}  \|e^x\phi\|_{L^2_x}\sum_{j=1}^\infty \alpha_j e^{4j}<+\infty.
\end{align*}

From \eqref{blowup_eq9} and \eqref{blowup_eq10} we conclude that
\begin{align*}
\sum_{j=1}^\infty \alpha_j \|e^{x} W(-j) \phi\|_{H^m(\mathbb R)}<+\infty.
\end{align*}

Now let us prove that, for every $t\in(0,+\infty)\setminus \mathbb Z^+$, $W(t)u_0\in C^1(\mathbb R)$. Initially, assume that $t>1$. Then there exists $j_0\in\mathbb Z^+$ such that $j_0<t<j_0+1$, and
\begin{align*}
W(t)u_0 = \sum_{j=1}^{j_0} \alpha_j W(t-j)\phi + \sum_{j=j_0+1}^\infty \alpha_j W(t-j) \phi.
\end{align*}

For $j\leq j_0$, $W(t-j)\phi = e^x \widetilde w (x,t-j)$, and thus, to see that $W(t-j)\phi\in C^1(\mathbb R)$, it is enough to prove that $\widetilde w(x,t-j)\in C^1(\mathbb R_x)$. This last statement is true if we prove that $\widetilde w(\cdot_x,t-j)\in H^2(\mathbb R_x)$. Indeed, for $j\leq j_0<t$, using \eqref{blowup_eq9}, we have that
\begin{align}
\notag \|\widetilde w(\cdot_x,t-j)\|_{H^2} \leq & C (\|\widetilde w(\cdot_x,t-j)\|_{L^2} + \|\partial_x^2 \widetilde w(\cdot_x,t-j)\|_{L^2})\\
\notag \leq & C (\| [\widetilde w(\cdot_x,t-j)]^\wedge(\cdot_\xi)\|_{L^2} + \|\xi^2 [ \widetilde w(\cdot_x,t-j)]^\wedge(\cdot_\xi)\|_{L^2})\\
\notag  \leq & C \left( \|e^{i(t-j)(-\xi^5+10\xi^3-5\xi)} e^{-5(t-j)(\xi^2-1)^2} e^{4(t-j)} (e^{-x}\phi)^\wedge(\xi)\|_{L^2_\xi} \right.\\
\notag& \left. + \|\xi^2 e^{i(t-j)(-\xi^5+10\xi^3-5\xi)} e^{-5(t-j)(\xi^2-1)^2}e^{4(t-j)}(e^{-x}\phi)^\wedge (\xi)\|_{L^2_\xi} \right)\\
\leq & C \left(e^{4(t-j)} \|e^{-x}\phi\|_{L^2_x} + e^{4(t-j)} \|\xi^2 e^{-5(t-j)(\xi^2-1)^2}\|_{L^\infty_\xi} \|e^{-x} \phi\|_{L^2_x} \right)<\infty.\label{blowup_eq11}
\end{align}

From \eqref{blowup_eq11}, we conclude that $\sum_{j=1}^{j_0} \alpha_j W(t-j)\phi \in C^1(\mathbb R)$. To show that
$$\sum_{j=j_0+1}^{\infty} \alpha_j W(t-j)\phi \in C^1(\mathbb R),$$
it suffices to demonstrate that
$$\sum_{j=j_0+1}^\infty \alpha_j e^x W(t-j) \phi \in C^1(\mathbb R).$$
For this, it is sufficient to prove that
$$\sum_{j=j_0+1}^\infty \alpha_j e^x W(t-j) \phi \in H^2(\mathbb R).$$
This last statement will be true if we prove that
\begin{align}
\sum_{j=j_0+1}^\infty \alpha_j \|e^x W(t-j)\phi\|_{H^2}<\infty. \label{blowup_eq12}
\end{align}

In fact,
\begin{align*}
&\sum_{j=j_0+1} \alpha_j \|e^x W(t-j) \phi\|_{H^2} = \sum_{j=j_0+1}^\infty \alpha_j \|w(\cdot_x,t-j)\|_{H^2}\\
& \leq C \sum_{j=j_0+1}^\infty \alpha_j (\|w(\cdot_x,t-j)\|_{L^2} + \|\partial_x^2 w(\cdot_x,t-j)\|_{L^2})\\
&\leq  C \sum_{j=j_0+1}^\infty \alpha_j \left( \|[w(\cdot,t-j)]^\wedge(\cdot_\xi)\|_{L^2} + \|\xi^2 [w(\cdot,t-j)]^\wedge(\cdot_\xi)\|_{L^2_\xi}\right)\\
&\leq  C \sum_{j=j_0+1}^\infty \alpha_j \left( e^{-4(t-j)}\|e^x \phi\|_{L^2_x} + \|\xi^2 e^{5(t-j)(\xi^2-1)^2} e^{-4(t-j)} (e^x \phi)^\wedge (\xi)\|_{L^2_\xi} \right)\\
&\leq  C e^{-4t} \sum_{j=j_0+1}^\infty \alpha_j \left( e^{4j}\|e^x \phi\|_{L^2_x} + e^{4j} \|\xi^2 e^{5(t-(j_0+1))(\xi^2-1)^2}  (e^x \phi)^\wedge (\xi)\|_{L^2_\xi} \right)\\
&\leq  C e^{-4t} \sum_{j=j_0+1}^\infty \alpha_j \left( e^{4j}\|e^x \phi\|_{L^2_x} + e^{4j} \|\xi^2 e^{5(t-(j_0+1))(\xi^2-1)^2}\|_{L^\infty_\xi} \|  e^x \phi \|_{L^2_x} \right)\\
&\leq  C e^{-4t} \|e^x \phi\|_{L^2_x} \sum_{j=j_0+1}^\infty \alpha_j  e^{4j}<+\infty,
\end{align*}
which is precisely \eqref{blowup_eq12}.\\

Let us observe that, if $0<t<1$, the analysis is the same as that carried out for the proof of the series $\sum_{j=j_0+1}^\infty \alpha_j W(t-j)\phi$.\\

Thus we conclude that if $t\in (0,+\infty)\setminus \mathbb Z^+$, it follows that
\begin{align*}
W(t)u_0 = \sum_{j=1}^{j_0} \alpha_j W(t-j)\phi + \sum_{j=j_0+1}^\infty \alpha_j W(t-j)\phi \in C^1(\mathbb R).
\end{align*}

Finally, let us see that, if $t=n$, for some $n\in \mathbb Z^+$, then $W(t)u_0$ belongs to $C^1(\mathbb R\setminus\{0\})\setminus C^1(\mathbb R)$. That is, let us see that
\begin{align*}
W(n)u_0 = \alpha_n \phi + \sum_{j\in \mathbb Z^+\setminus \{n\}} \alpha_j W(n-j) \phi \equiv \alpha_n \phi + \Phi_n \in C^1(\mathbb R\setminus \{0\}) \setminus C^1(\mathbb R).
\end{align*}

Proceeding as in the case $t\notin \mathbb Z^+$, it can be shown that
\begin{align*}
\Phi_n \equiv \sum_{j\in \mathbb Z^+ \setminus \{n\}} \alpha_n W(n-j) \phi\in C^1(\mathbb R),
\end{align*}

and since $\alpha_n \phi \in C^1(\mathbb R\setminus \{0\})\setminus C^1(\mathbb R)$, then $W(n)u_0 \notin C^1(\mathbb R)$. Otherwise, since $\alpha_n>0$, we would have that $\phi\in C^1(\mathbb R)$. Moreover, since $\alpha_n \phi \in C^1(\mathbb R\setminus \{0\})$, and $\Phi_n \in C^1(\mathbb R\setminus \{0\})$, it follows that $W(n)u_0 \in C^1(\mathbb R\setminus \{0\})$.\\

Theorem \ref{blowup} is proved.\qed\\

\textbf{Acknowledgments}\\

This work was partially supported by Universidad Nacional de Colombia, Sede-Medellín– Facultad de Ciencias – Departamento de Matemáticas – Grupo de investigación en Matemáticas de la Universidad Nacional de Colombia Sede Medellín, carrera 65 No. 59A -110, post 50034, Medellín Colombia. Proyecto: Análisis no lineal aplicado a problemas mixtos en ecuaciones diferenciales parciales, código Hermes 60827. Fondo de Investigación de la Facultad de Ciencias empresa 3062.

\end{document}